\documentclass[a4paper]{article}

\usepackage{amsmath,amssymb,amsthm, amscd}
\usepackage{array}
\usepackage{dsfont}
\usepackage{epsfig}

\usepackage{tikz}

\begin{document}

\newcommand{\RT}{Reshetikhin-Turaev }
\newcommand{\wrt}{Witten-Reshetikhin-Turaev }
\newcommand{\Mod}{\operatorname{Mod}}
\newcommand{\Emcg}{\widetilde{\Mod(\Sigma)}}
\newcommand{\lollipopgraph}{
\tikz[scale=0.5,baseline=-0.5ex]{
\draw (0,0) circle (1);
\draw (-3,0)--(-1,0);
\draw (-2,0) node[above]{2c};
\draw (0,1) node[above]{i+c};
}
}
\newcommand{\U}{\operatorname{U}}
\newcommand{\quotient}[2]{{\raisebox{.2em}{$#1$}\left/\raisebox{-.2em}{$#2$}\right.}}
\newcommand{\GL}{\operatorname{GL}}
\newcommand{\Thetagraph}{
\tikz[scale=0.3,baseline=-0.5ex]{
\draw (0,0) circle (2);
\draw (0,2)--(0,-2);
}
}
\newcommand{\emcg}{\widetilde{\Mod(\Sigma_g)}}

\newcommand{\Twicegraph}{
\tikz[scale=0.5,baseline=-0.5ex]{
\draw (0,0) circle (1);
\draw (-3,0)--(-1,0);
\draw (1,0)--(3,0);
\draw (-2,0) node[above]{$1$};
\draw (0,1) node[above]{$i+\frac{c-1}{2}$};
\draw (0,-1) node[below]{$i+\frac{c+1}{2}$};
\draw (2,0) node[above]{$c$};
}
}

\pagestyle{plain}

\theoremstyle{plain}
\newtheorem{theorem}{Theorem}[section]
\newtheorem{main_theorem}[theorem]{Main Theorem}
\newtheorem{proposition}[theorem]{Proposition}
\newtheorem{corollary}[theorem]{Corollary}
\newtheorem{corollaire}[theorem]{Corollaire}
\newtheorem{lemma}[theorem]{Lemma}
\theoremstyle{definition}
\newtheorem{definition}[theorem]{Definition}
\newtheorem{Theorem-Definition}[theorem]{Theorem-Definition}
\theoremstyle{remark}
\newtheorem*{remark}{Remark}
\newtheorem*{Remark}{Remark}
\newtheorem{example}{Example}

\sloppy

\title{On the (in)finiteness of the image of Reshetikhin-Turaev representations}

\author{Julien \textsc{Korinman}
\\ \small Universidade Federal de S\~ao Carlos
\\ \small Rodovia Washington Lu\`is, Km 235, s/n 
\\ \small S\~ao Carlos - SP, 13565-905
\\ \small email: \texttt{ julienkorinman@dm.ufscar.br}
}
\date{}
\maketitle


\begin{abstract} 
We state a simple criterion to prove the infiniteness of the image of \RT irreducible representations of the mapping class groups of surfaces. We use it to study some of the \RT representations associated to the tori with one and two punctures and derive an alternative proof of the results of \cite{Fu99}. 
\vspace{2mm}
\par
Keywords: \wrt representations, mapping class group,  Topological Quantum Field Theory.
\vspace{2mm}
\par  AMS 2010 Subject Classification: 57M99
\end{abstract}

\section{Introduction}
\vspace{3mm}
\par  Witten  constructed in \cite{Wi2} a $2+1$ dimensional TQFT using path integrals  which gives a three-dimensional interpretation of the Jones polynomial and produces invariants for $3$-dimensional closed oriented manifolds equipped with  framed links and additional structure. Reshetikhin and Turaev (\cite{RT}) gave a rigorous construction using quantum groups and latter Blanchet, Masbaum, Habegger and Vogel constructed these TQFTs by means of  Kauffman bracket skein algebra (\cite{BHMV2}) following the work of Lickorish (\cite{Li2}). We will use the construction in \cite{BHMV2}.
\vspace{3mm}
\par These TQFTs give rise to finite dimensional representations $\rho_p$, $p\geq 3$ of some central extension of the mapping class group $\Emcg$ of any closed oriented surface $\Sigma$,  equipped with colored points.  Here we will discuss whether these representations have finite image or not. Several studies have been made in that sense.
The \RT representations associated to a torus without marked points have finite image. This fact was known in the conformal field theory community (see \cite{CG} and references herein) and has been proved independently by Gilmer in \cite{Gi99}. In higher genus, the Reshetikhin-Turaev representations have finite image at level $3$ and $6$ (see \cite{Wr96}) and they are characters at level $4$ . In any other cases, they have infinite image. It was proved by Funar in \cite{Fu99} when $(g,p)\neq (2, 20)$ (see also \cite{Ma99} for a proof that they contain an element of infinite order). The last remaining case $g=2$, $p=20$ was treated in \cite{EF}. Concerning the representation associated to a one holed torus, it results from \cite{Sa12} that they have infinite image for high enough level when the level is odd and the dimension is fixed.

\vspace{3mm}
\par In this paper, we develop a simple criterion to check whether a given Reshetikhin-Turaev representation has infinite image or not when this representation is irreducible. It permits us to recover the above results in these cases, and to study the finiteness of the representations associated to a one-holed torus, which could not be derived from previous papers.  
\vspace{2mm}
\par Denote by $\rho_p^c$ the level $p$ Reshetikhin-Turaev representation associated to a torus $\mathcal{T}^c$ equipped with a puncture colored by $2c$ and by $\rho_p^{1,c}$ the representation associated to a torus $\mathcal{T}^{(1,c)}$ with two punctures colored by $1$ and $c$ respectively. The main theorem is the following:

\begin{theorem}\label{th_infinite_torus} Let $r\geq 4$. Then we have:
\begin{enumerate}
\item If $2c=r-2$ or $2c=r-3$ and $r$ is odd, then $\rho_{2r}^c$ has finite image.
\item If $2c<r-3$ and $r$ is an odd prime, then $\rho_{2r}^c$ has infinite image.
\item If $r>5$ is prime and  $r\equiv 3 \pmod{8}$ or $r\equiv 5 \pmod{8}$, then $\rho_{r}^c$ has infinite image.
\item If $r\geq 5$, then $\rho_{2r}^{(1,1)}$ has infinite image.
\item If $r\geq 5$ is odd and $c\equiv 1 \pmod{4}$, then $\rho_{2r}^{(1,c)}$ has infinite image.
\end{enumerate}
\end{theorem}

\vspace{3mm}
 \textbf{Acknowledgements:} The author is thankful to  L.Funar and F.Costantino for useful discussions and to the referee for important remarks. He acknowledges  support from the grant ANR $2011$ BS $01 020 01$ ModGroup and the NSF grants DMS-$1107452$,$ 1107263$ and $1107367$ RNMS: GEometric structures And Representation varieties (the GEAR Network).

\vspace{3mm}
\section{Reshetikhin-Turaev representations and basis of conformal blocks}
\vspace{3mm}
\par Let $\Sigma$ be a closed oriented surface with colored punctures and $p\geq 3$. The associated representation $\rho_p$ acts on a finite rank $k_p$-module $V_p(\Sigma)$ where $k_p$ is the quotient of $\mathbb{Z}[A,\frac{1}{p},\kappa]$  by the relations $\phi_{2p}(A)=0$, where $\phi_{2p}$ represents the $2p$-th cyclotomic polynomial, and $\kappa^6=A^{-6-\frac{p(p+1)}{2}}$. The module $V_p(\Sigma)$ is equipped with a non-degenerate Hermitian form $\left<\cdot, \cdot\right>_p$ valued in $k_p$ which is preserved by the action of $\Emcg$. Note that in \cite{BHMV2}, the invariant pairing is bilinear whereas in this paper we choose it Hermitian using the involution of $k_p$ sending $A$ to $A^{-1}$. 

\vspace{3mm}
\par  Note $P$ the set of punctures in $\Sigma$. Let $(\gamma_e)_e$ be a pants decomposition of $\Sigma \backslash P$, that is a maximal set of isotopy classes of non intersecting non-contractible pairwise non-homotopic simple closed curves in $\Sigma \backslash P$. Such a pants decomposition can be equivalently described by a banded uni-trivalent graph $\Gamma$ embedded in a handlebody bounded by $\Sigma$ (see \cite{BHMV2} for details). We denote $E(\Gamma)$ and $V(\Gamma)$ the sets of edges and trivalent vertices respectively of $\Gamma$. We associate to any such $\Gamma$ a basis $(u_{\sigma})_{\sigma}$ of $V_p(\Sigma)$ which is orthogonal for the invariant form $\left<\cdot, \cdot\right>_p$. Moreover the vectors $u_{\sigma}$ are common eigenvectors for the image of (any lift in $\Emcg$ of) the Dehn twists $\rho_p(T_{\gamma_e})$ associated to the pants decomposition. Let $I_p$ denotes the set $\{0, \ldots, r-2\}$ if $p=2r$, and the set of even integers of $\{0, \ldots, r-2\}$ if $p=r$ is odd. The vectors $u_{\sigma}$ of the basis of conformal blocks are indexed by $p$-admissible colorings of $\Gamma$, that is by maps $\sigma : E(\Gamma)\rightarrow I_p$ such that: 
\begin{itemize}
\item If $v$ is a puncture of $\Sigma$ colored by $c_v$ and $e_v\in E(\Gamma)$ is the corresponding edge, then $\sigma(e_v)=c_v$.
\item If $e_1, e_2, e_3 \in E(\Gamma)$ are three edges adjacent to a common vertex, then $\sigma(e_1)+\sigma(e_2)+\sigma(e_3)$ is an even integer strictly  smaller than $2r-2$.
\item If $e_1, e_2, e_3 \in E(\Gamma)$ are three edges adjacent to a common vertex, then $ \sigma(e_1)\leq \sigma(e_2)+\sigma(e_3)$.
\end{itemize}
\vspace{3mm}
\par Define the quantum numbers by $[n]:=\frac{A^{2n}-A^{-2n}}{A^2-A^{-2}}\in k_p$ and $[n]!:=[n][n-1]\ldots [1]\in k_p$. Let $\sigma$ be a $p$-admissible coloring of $\Gamma$. If $e$ is an edge of $\Gamma$, we note $<\sigma(e)>:=(-1)^{\sigma(e)}[\sigma(e)+1]\in k_p$. If $v$ is a trivalent vertex of $\Gamma$ with adjacent edges $e_1, e_2, e_3$, set $i:=\frac{\sigma(e_1)+\sigma(e_2)-\sigma(e_3)}{2}$, $j:=\frac{\sigma(e_1)-\sigma(e_2)+\sigma(e_3)}{2}$ and $k:=\frac{-\sigma(e_1)+\sigma(e_2)+\sigma(e_3)}{2}$. We then define:
$$ <\sigma(v)>:= (-1)^{i+j+k} \frac{[i+j+k+1]![i]![j]![k]!}{ [\sigma(e_1)]![\sigma(e_2)]![\sigma(e_3)]!}$$

\vspace{2mm}
\par An important property we will use in this paper is the following equality (Theorem $4.11$ in \cite{BHMV2}):
\begin{equation}\label{equation1}
\left< u_{\sigma}, u_{\sigma}\right> = \eta^{\# V(\Gamma) -\#E(\Gamma)} \frac{\prod_{v\in V(\Gamma)}<\sigma(v)>}{\prod_{e\in E(\Gamma)} <\sigma(e)>}
\end{equation}
\vspace{2mm}
where  $\eta=\frac{1}{2p}(A\kappa)^3(A^2-A^{-2})\sum_{m=1}^{2p}{(-1)^mA^{-m^2}}\in k_p$ is an invertible element which represents the sphere invariant in TQFT. Note that the right hand side is a product of  quantum numbers $[n]$ with $n\in I_p$ and their inverses and that changing $\kappa$ to $-\kappa$ changes the signature of $\left<\cdot,\cdot\right>$ to its opposite.
\vspace{3mm}
\par A choice of a particular $2p$-th primitive root of unity $A\in \mathbb{C}$ and a compatible complex $\kappa\in\mathbb{C}$ gives a complex vector space $V_p^{A,\kappa}(\Sigma)=V_p(\Sigma)\otimes\mathbb{C}$ and a complex representation $\rho_p^{A,\kappa}$ which preserves the non-degenerate Hermitian form $\left<\cdot,\cdot \right>_p^{A,\kappa}$. 
Note that if $A=\exp\left(\frac{i\pi l}{p}\right)$ with $l$ prime to $2p$, then the quantum numbers become $[n]=\frac{\sin\left(\frac{2i\pi l n}{p}\right)}{\sin\left(\frac{2i\pi l}{p}\right)}\in \mathbb{R}$. In particular they are nonzero if and only if  $n$ is not a multiple of $p$ and equation $\eqref{equation1}$ shows the non-degeneracy of the invariant form. Moreover if $A=\exp\left(\frac{i\pi }{p}\right)\in \mathbb{C}$ and $p$ is even, then for any $0\neq n\in I_p$, we have $[n]=\frac{\sin\left(\frac{2i\pi  n}{p}\right)}{\sin\left(\frac{2i\pi }{p}\right)}>0$. Formula \ref{equation1} implies that $\left<\cdot, \cdot\right>_p^{A,\kappa}$ is positive definite when $\kappa$ is such that $\eta>0$. This can be used to prove the classical fact that that the \RT representations are semi-simple.

\vspace{3mm}
\par In the particular case where $\Sigma=\mathcal{T}^c$ is a torus with one puncture colored by $2c\in I_p$ and $\Gamma$ is the standard Lollipop graph, $V_p(\mathcal{T}^c)$ has an orthogonal basis $(u_i^c)_i$ associated to the colorings sending the stick edge to $2c$ and the loop edge to $i+c$, that is:
$$ u_i^c := \lollipopgraph $$
\par Here the index $i$ lies in the set $\{0, \ldots, r-2-2c\}$ if $p=2r$ and in the subset of even elements of $\{0, \ldots, r-2-2c\}$ if $p=r$ is odd. Formula \eqref{equation1} reads:
\begin{equation}\label{equation2}
\left<u_i^c, u_i^c\right> = \eta^{-1}(-1)^c\frac{[2c+i+1]![i]! ([c]!)^2}{[i+c+1]![i+c]![2c+1]!}
\end{equation}
\par
\begin{figure}[!h] 
\centerline{\includegraphics[width=8cm]{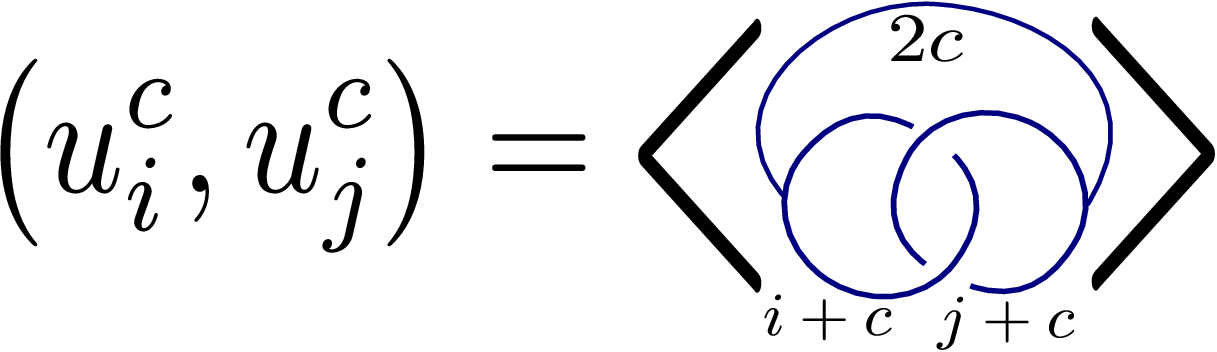} }
\caption{Definition of the Hopf pairing for the one-holed torus.} 
\label{fig_HopfPairingTorus} 
\end{figure}

 In \cite{GM11}, the authors defined another Hermitian form on $V_p(\mathcal{T}^c)$, named Hopf pairing, and denoted $(\cdot, \cdot)_p$. This form is defined by the property that $\left( u_i^c, u_j^c\right)_p$ is the colored banded graph Kauffman-bracket invariant of the graph drawn in Figure \ref{fig_HopfPairingTorus}. Note that it is not invariant under the action of $\widetilde{\Mod}(\mathcal{T}^c)$.

\vspace{2mm}
\par We will use the important fact that there exists an element $\phi\in \widetilde{\Mod}(\mathcal{T}^c)$ such that $(x,y)_p = \left<x, \rho_p^c(\phi)y\right>_p$ for all $x,y\in V_p(\mathcal{T}^c)$. Moreover we will use the following:
\begin{lemma}\label{lemma_Masbaum}
For any $ i\in \{0, \ldots, r-2-2c\}$, we have $\left(u_0^c, u_i^c\right)_p \neq 0$.
\end{lemma}
\vspace{3mm}
\begin{proof}
In the proof of Theorem $2.2$ of \cite{GM11}, the authors computed:
$$ (u_0^c, u_i^c) = C_{c+i, c, c} \left< \begin{array}{lll} 2c & c+i & c+i \\ i & c & c \end{array} \right> $$
 where we have the formula from the proof of Lemma $4.1$ in \cite{Masbaum03} p.$550$:
$$C_{c+i, c, c}=A^{ci}\prod_{j=i+1, \ldots, i+c}(1-A^{-2j})$$ 
  \vspace{3mm}
    and  $\left< \begin{array}{lll} 2c & c+i & c+i \\ i & c & c \end{array} \right>$ is a tetrahedron coefficient whose formula is given in Theorem $2$ of \cite{MV} and is:
  $$ \left< \begin{array}{lll} 2c & c+i & c+i \\ i & c & c \end{array} \right> = (-1)^i\frac{([c]!)^2 [i]! [2c+i+1]!}{[2c]!([c+i]!)^2} $$ 
   \par Both these numbers are product of  quantum numbers  $ [n]=\frac{A^{2n}-A^{-2n}}{A^2-A^{-2}}\in k_p$ (and their inverse) with $1\leq n \leq r-2$ and in $(1-A^{2j})$ for $1\leq j\leq r-2$. In particular they are not null. 
\end{proof}

\vspace{3mm}
\section{Complete positivity}
\vspace{3mm}

 \begin{definition}
We say that $V_p(\Sigma)$ is \textit{completely positive} if for any choice of $A\in \mathbb{C}$ as a primitive $2p$-th root of unity and compatible $\kappa\in \mathbb{C}$, the Hermitian form $(V_p^{A,\kappa}(\Sigma),\left<\cdot,\cdot \right>_p^{A,\kappa})$ is positive definite or negative definite.
\end{definition}
\vspace{3mm}
\par  The only effect of changing $\kappa$ is possibly to change the eigenvalues of the Hermitian form to its opposite, by changing the sign of $\eta$ in Formula \eqref{equation1}, so complete positivity only means that all eigenvalues have the same sign for a given $\kappa$ and all $A$. 
\vspace{4mm}
\begin{proposition}\label{prop_CP}
\begin{enumerate}
\item
Let $r$ be an odd prime number and $p=r$ or $p=2r$. Then $(V_p(\Sigma),\left<\cdot,\cdot \right>_p)$ is completely positive if and only if $\rho_p$ has finite image.
\item  If  $(V_p(\mathcal{T}^{(1,c)}), \left<\cdot,\cdot \right>_p)$ is not completely positive, then $\rho_p^{(1,c)}$ has infinite image.
\end{enumerate}
\end{proposition}
\vspace{3mm}

\par We first show the second part of Proposition \ref{prop_CP}. The proof relies on the following lemma made by Coxeter in \cite{Coxeter} p.$116$ (see also \cite{Coxeter_Moser} p.$121$):
\begin{lemma}[Coxeter \cite{Coxeter} p.$116$]\label{lemma_Coxeter}
\par Let $V$ be a finite complex vector space equipped with a non-degenerate Hermitian form $\left<\cdot,\cdot\right>$ which is neither negative definite nor positive definite. Let $G\subset\U(V, \left<\cdot,\cdot\right>)$ be a group acting linearly by preserving the Hermitian form. Suppose $G$ acts irreducibly on $V$, then $G$ is infinite.
\end{lemma}
\vspace{3mm}
\par Since Coxeter's argument is short and stated slightly differently, we briefly explain it. Let $\left(\cdot, \cdot\right)$ be any Hermitian positive definite form on $V$. Suppose by contradiction that $G$ is finite. We define a $G$-invariant positive definite Hermitian form on $V$ by the formula $\left< u,v\right>':= \frac{1}{\# G} \sum_{g\in G}\left( g\cdot u, g\cdot v\right)$. To any $t\geq 0$ we associate the $G$-invariant form $\left<\cdot,\cdot\right>_t := \left<\cdot,\cdot\right>' +t \left<\cdot,\cdot\right>$. Recall the classical result (see e.g. \cite{Franklin} Theorem $2$ p.$257$) that we can always co-diagonalize  two Hermitian forms when one is positive definite. Let $(e_i)_i$ be a basis of $V$ which is orthonormal for $\left<\cdot,\cdot\right>'$ and orthogonal for $\left<\cdot,\cdot\right>$ and write $\lambda_i:=\left<e_i, e_i\right>$. Since $\left<\cdot,\cdot\right>$ is indefinite, some $\lambda_i$ are positive and some are negative. Choose $i_0$  such that $\lambda_{i_0}<0$ is minimal among the $\lambda_i$ and note $t_0:=-\frac{1}{\lambda_{i_0}}$. Then $\left<e_i, e_i\right>_{t_0} = 1- \frac{\lambda_i}{\lambda_{i_0}}$ is non-negative and null  if and only if $\lambda_i=\lambda_{i_0}$. Thus the kernel of the $G$-invariant form $\left<\cdot,\cdot\right>_{t_0}$ is $G$-invariant, non-trivial and proper. This contradicts the assumption of irreducibility of the action of $G$ on $V$.

\vspace{3mm}
\begin{proof}[Proof of the second part of Proposition \ref{prop_CP}]

\par Suppose that $(V_p(\Sigma), \left<\cdot,\cdot \right>)$ is not completely positive and fix $A\in\mathbb{C}$ and a compatible $\kappa\in\mathbb{C}$ such that $\left<\cdot, \cdot\right>^{A,\kappa}$ is nor positive nor negative definite. If the associated representation is irreducible, then Lemma \ref{lemma_Coxeter} implies that the representations has infinite image. It follows from Theorem $1.1$ in  \cite{KoberdaSantharoubane17} that the representations $\rho_p^{(1,c)}$ are irreducible for every levels $p$.  

 When $p$ is an odd prime, Corollary $3.2$ of \cite{GM12} implies that the representation associated to any punctured surface is irreducible. When $p=2r$ with $r$ an odd prime, in \cite{BHMV2} the authors showed that there exists a tensor decomposition $\rho_{2r}^c\cong \rho_r^c\otimes \rho'_2$ where $\rho'_2$ is irreducible and totally positive. Hence $\rho_r$ is also totally positive and both $\rho_r$ and $\rho_{2r}$ have infinite image.

 However in the one-holed torus case,  a more elementary proof of the irreducibility can be derived by adapting Roberts' argument of \cite{Ro}. Indeed consider the Lollipop basis $(u_i^c)_i$ described in the previous section. The assumption that $p=r$ or $p=2r$ with $r$ an odd prime implies that the operator associated by $\rho_p^c$ to a Dehn twist along the meridian has eigenvalues with multiplicity one. Moreover the $u_i^c$ are eigenvectors of this operator. As a consequence, every irreducible subspace of $V_p(\mathcal{T}^c)$ is spanned by vectors of the Lollipop basis. Now Lemma \ref{lemma_Masbaum} implies that $u_0^c$ is cyclic. Together with the fact that the eigenvalues have multiplicity one, this proves the irreducibility of $\rho_p^c$.

\end{proof}

\vspace{6mm}
\par When $p$ is an odd prime or   $p=2r$ with $r$ an odd prime,
\\ we write $\alpha_p:= \left\{ \begin{array}{l} 2p  \mbox{, if }p\equiv 3\pmod{4} \\ 4p  \mbox{, if }p\equiv 1,2 \pmod{4} \end{array}\right.$ and $\mathcal{O}_p:=\quotient{\mathbb{Z}[A]}{\phi_{\alpha_p}(A)}$.
\vspace{1mm}
\par
 It was showed in \cite{GM07, Qa13} that $V_p(\Sigma)$ contains a free $\mathcal{O}_p$ lattice of maximal rank preserved by $\Emcg$.
\vspace{3mm}
\par We denote by $\mu(\alpha_p)=\{ q_1, \ldots, q_{\varphi(\alpha_p)} \}$ the set of primitive $\alpha_p$-th roots of unity and by $\varphi(\alpha_p)=\# \mu(\alpha_p)$ the Euler totient function. The canonical embedding is the injective linear map: 
$$\Psi : \mathcal{O}_p \rightarrow \mathbb{C}^{\varphi(\alpha_p)} $$
 sending $A^n$ to $\left(q_1^n, \ldots, q_{\varphi(\alpha_p)}^n\right)$. It is well known that its image is a discrete lattice.

\vspace{4mm} \par 
\begin{proof}[Proof of the first part of Proposition \ref{prop_CP}]
\par We suppose that $p=r$ or $p=2r$ with $r$ an odd prime.
\par A classical result of Roberts (\cite{Ro}) states that the representations $\rho_p$ are irreducible. Coxeter's lemma \ref{lemma_Coxeter} thus implies that if the space is not totally positive, then the representation has infinite image.

\vspace{3mm}
\par Conversely, if $V_p(\Sigma)$ is completely positive, using the $\mathcal{O}_p$ lattice of \cite{GM07,Qa13}, we have an injective homomorphism from $\rho_p(\Emcg)$ to the group of matrices with coefficients in $\mathcal{O}_p$. Once composed with the canonical embedding  $\Psi$, we get an injective group morphism $$\widetilde{\Psi} : \rho(\Emcg)\rightarrow \GL_d(\mathbb{C})\times \ldots \times \GL_d(\mathbb{C})$$ 
 where $d$ denotes the dimension of $V_p(\Sigma)$. Since the image of the canonical embedding is discrete, so is the image of  $\widetilde{\Psi}$. Moreover the complete positivity of $V_p(\Sigma)$ implies that the image lies in the compact product of $r-1$ unitary groups. This implies the finiteness of $\rho_p(\Emcg)$.
\end{proof}

\vspace{5mm}
\section{(In)finiteness of the holed torus representations}
\vspace{3mm}
Using Proposition $\ref{prop_CP}$, we can prove the main theorem of this paper:
\vspace{3mm}
\begin{proof}[\textit{(Proof of Theorem \ref{th_infinite_torus})}]
\par We first consider the torus with one puncture. When $p=2r$ and $2c=r-2$, the representation $\rho_p^c$ is a character. The group  $ \Mod(\mathcal{T}^c)$ is generated by two Dehn twists $T$ and $T'$ along a longitude and a meridian of the torus. Since the images of these two elements have finite order, the representation has finite image.
\vspace{2mm}
\par Next when $p=2r$ and $2c=r-3$ we could show that the $2$-dimensional spaces $V_p(\mathcal{T}^c)$ are completely positive but Proposition \ref{prop_CP} would permit us to conclude only when $r$ is prime. Instead we use the following argument. First since the image is semi-simple and non-abelian, the representation is irreducible. In \cite{Fo96},  Formanek proved  that the only $2$-dimensional irreducible representations of $B_3=\left< t,t' | tt't=t'tt'\right>$ are conjugate to one of the $\chi(y)\otimes\beta_3(q)$, where $\chi(y)$ is the character which sends both generators $t$ and $t'$ to $y\in\mathbb{C}^*$ and $\beta_3(q)$ is the reduced Burau representation in $q$. Since the eigenvalues of $\rho_p^c(t)$ are $\mu_c=(-1)^cA^{c(c+2)}$ and $\mu_{c+1}=(-1)^{c+1}A^{(c+1)(c+3)}$ we must have $y=-\mu_c$ and $q=-A^r$ or $y=-\mu_{c+1}$ and $q=-A^{-r}$ so $-q$ is a $4-{th}$ primitive root of unity. Moreover such Burau representations $\beta_3(-q)$ at $4-th$ roots of unity was proven to have finite image in \cite{FK14} Proposition $3.1$.

\vspace{3mm}
\par Now if $p=2r$ and $2c<r-3$, using Proposition $\ref{prop_CP}$, to show that $\rho_p^c$ has infinite image we need to find a $2p-{th}$ root of unity $A$ and $1\leq i \leq r-3-2c$ such that $\frac{\left<u_i^c, u_i^c\right>}{\left<u_0^c, u_0^c\right>}(A)$ is negative. Write $A=\exp\left(\frac{i\pi k}{p}\right)$ with $g.c.d.(k, 2p)=1$.
\vspace{2mm}
\par If $p=2r$ with $r$ odd and $2c<r-3$, using Equation $\eqref{equation2}$, we compute:
\begin{eqnarray*}
\frac{\left<u_2^c, u_2^c\right>}{\left<u_0^c, u_0^c\right>} &=& \frac{[2c+3][2c+2][2]}{[c+1][c+3][c+2]^2} 
\\ &=& \frac{\sin\left(\frac{\pi k}{r}(2c+3)\right) \sin\left(\frac{\pi k}{r}(2c+2)\right) \cos\left(\frac{\pi k}{r}\right) \sin\left(\frac{\pi k}{r}\right)^2}{\sin\left(\frac{\pi k}{r}(c+1)\right) \sin\left(\frac{\pi k}{r}(c+3)\right) \sin\left(\frac{\pi k}{r}(c+2)\right)^2}
\end{eqnarray*} 
\par Choosing $k=\left\{ \begin{array}{ll} \frac{r-1}{2} & \mbox{, if } r\equiv 3 \pmod{4} \\ \frac{r+1}{2} & \mbox{, if } r\equiv 1 \pmod{4}\end{array} \right. $, we have $\frac{\left<u_2^c, u_2^c\right>}{\left<u_0^c, u_0^c\right>} <0$.

\vspace{2mm}
\par Finally when $p\equiv 3\pmod{8}$ choose $k=\frac{p+1}{4}$, when $p\equiv 5\pmod{8}$ choose $k=\frac{p-1}{4}$. Then $\frac{\left<u_2^c, u_2^c\right>}{\left<u_0^c, u_0^c\right>}$ has the same expression than in the previous case (replace $r$ by $p$) and is negative for the same reason.
\vspace{2mm}
\par We now consider the torus with two punctures colored by $1$ and $c$ with $c$ odd and write $p=2r$ with $r\geq 5$ and $A=\exp\left(\frac{i\pi k}{p}\right)$. Given $i\in \{0, \ldots, r-\frac{c+5}{2}\}$, define the vector 
$u_i^{(1,c)}:= 
\Twicegraph 
\in V_p(\mathcal{T}^{(1,c)})$.
 We compute:
$$ \frac{\left<u_1^{(1,c)}, u_1^{(1,c)}\right>}{\left<u_0^{(1,c)}, u_0^{(1,c)}\right>} = \frac{[c+2][3]}{[\frac{c+3}{2}][\frac{c+5}{2}][2]}
 = \frac{\sin\left(\frac{\pi k}{r}(c+2)\right)\sin\left(\frac{3\pi k}{r}\right)\sin\left(\frac{\pi k}{r}\right)}{\sin\left(\frac{\pi k(c+3)}{2r}\right)\sin\left(\frac{\pi k(c+5)}{2r}\right)\sin\left(\frac{2\pi k}{r}\right)}$$
 If $c=1$, this expression has the sign of $[3]$, which is negative for $k=r-1$ when $r$ is even and for $k=r-2$ when $r$ is odd. If $c\equiv 1 \pmod{4}$ and $r$ is odd, then setting $k=\left\{ \begin{array}{ll} \frac{r-1}{2} & \mbox{, if } r\equiv 3 \pmod{4} \\ \frac{r+1}{2} & \mbox{, if } r\equiv 1 \pmod{4}\end{array} \right. $ we see that the above expression is negative and we conclude likewise.

\end{proof}
\vspace{3mm}

\begin{remark} When $p=2r$ and  $2c=r-4$ (i.e. when the representation is $3$-dimensional), then $V_p(\mathcal{T}^c)$ is completely positive. Indeed, putting $A=\exp \left(\frac{i\pi k}{r}\right)$ with $g.c.d.(k,2r)=1$, we have:
$$ \frac{\left<u_1^c, u_1^c\right>}{\left<u_0^c, u_0^c\right>}= \frac{\left<u_2^c, u_2^c\right>}{\left<u_1^c, u_1^c\right>} = \frac{ \sin\left( \frac{2\pi k}{r} \right)\sin\left( \frac{\pi k}{r} \right)}{\cos\left( \frac{\pi k}{r} \right)}> 0$$

\par We think the image should be finite in this case.
\end{remark}

\vspace{5mm}
\subsection{(In)finiteness of Reshetikhin-Turaev representations associated to closed surfaces}
\vspace{3mm}
\par The following theorem results from \cite{Fu99, EF}. We derive another proof from Proposition \ref{prop_CP}. Here we denote by $\rho_{p,g}$ the \RT representation at level $p$ associated to a genus $g$ closed surface $\Sigma_g$ without punctures.

\begin{theorem}[\textit{(\cite{Fu99, EF})}] \par 
Let $g\geq 2$, $r\geq 5$ (not necessary prime),  then $\rho_{2r,g}$ has infinite image.
\end{theorem}
\begin{proof}

\par Write $\Sigma_g \cong \Sigma_1^2 \bigcup_{S^1\cup S^1}\Sigma_{g-2}^2$. Then using Theorem $1.14$ of \cite{BHMV2}, we have:
$$ V_p(\Sigma_g)\cong \oplus_{c_1,c_2\in \{0, \ldots, r-2\} } V_p(\mathcal{T}^{(c_1,c_2)})\otimes V_p(\Sigma_{g-2}^{(c_1,c_2)}) $$

\par Any homeomorphism of $\Sigma_1^2$ uniquely extends to a homeomorphism of $\Sigma_g$ which is the identity on $\Sigma_{g-2}^2$. This gives an embedding $\widetilde{\Mod}(\mathcal{T}^{(c_1,c_2)})\rightarrow \emcg$ and the image of $\rho_p^{(1,1)}$ embeds in the image of $\rho_{p,g}$. Theorem $\ref{th_infinite_torus}$ states that this image is infinite.

\end{proof}


\bibliographystyle{alpha}
\bibliography{biblio}

\end{document}